\documentclass[10pt]{article}
\usepackage{amsmath,amsthm,amsfonts}

\newtheorem{theorem}{Theorem}[section]
\newtheorem{lemma}[theorem]{Lemma}

\newcommand{\definition}[1]{{\it #1}}
\newcommand{\dual}{^{\vee}}
\newcommand{\ses}[3]{
0 \rightarrow #1 \rightarrow #2 \rightarrow #3 \rightarrow 0
}

\DeclareMathOperator{\ddim}{\textbf{dim}}
\DeclareMathOperator{\ext}{ext}
\DeclareMathOperator{\Hom}{Hom}
\DeclareMathOperator{\Ext}{Ext}
\DeclareMathOperator{\End}{End}

\DeclareMathOperator{\hcf}{hcf}
\DeclareMathOperator{\E}{\textbf{E}}
\DeclareMathOperator{\Rep}{\textbf{Rep}}
\DeclareMathOperator{\Coh}{\textbf{Coh}}

\begin{document}

\title{Birational classification of moduli spaces of vector bundles
over $ \mathbb{P}^{2}$ }
\author{Aidan Schof\i eld}
\maketitle

\begin{abstract}
The depth of a vector bundle $E$ over $\mathbb{P}^{2}$ is the largest
integer $h$ such that $[E]/h$ is in the Grothendieck group of coherent
sheaves on $\mathbb{P}^{2}$ where $[E]$ is the class of $E$ in this
Grothendieck group. We show that a moduli space of vector bundles is
birational to a suitable number of $h$ by $h$ matrices up to
simultaneous conjugacy where $h$ is the depth of the vector bundles
classified by the moduli space. In particular, such a moduli space is
a rational variety if $h\leq 4$ and is stably rational when $h$
divides $420$.  
\end{abstract}

\section{Introduction}
\label{Intro}

The purpose of this paper is to study moduli spaces of vector bundles
over $\mathbb{P}^{2}$ birationally. Particular cases are known to be
rational and other cases, for example, the moduli space of vector
bundles of rank $n$, first Chern class $0$ and second Chern class $n$
are known to be birational to two $n$ by $n$ matrices up to
simultaneous conjugacy. We shall see that this is a general
phenomenon, that is, any such moduli space is birational to a suitable
number of suitably sized matrices up to simultaneous conjugacy. Our
method will be to reduce to a problem for representations of a
suitable quiver with relations and then to apply the results and
methods of \cite{qmod}. The results of this paper depend heavily on
those of \cite{qmod} and the reader will be assumed to have some
familarity with this paper.

From the work of Beilinson \cite{Beilinson}, one knows that the
category of vector bundles over $\mathbb{P}^{2}$ is derived equivalent
to the category of representations of a quiver with relations; this is
the quiver with $3$ vertices $u$, $v$ and $w$ and $3$ arrows from $u$
to $v$, which are $x$, $y$ and $z$, and $3$ arrows from $v$ to $w$
which are $x'$, $y'$ and $z'$ with relations $xy'=yx'$, $xz'=zx'$ and
$yz'=zy'$. The path algebra of this quiver with relations is
the endomorphism ring of $\mathcal{O}\oplus \mathcal{O}(-1)\oplus
\mathcal{O}(-2)$. Thus, roughly speaking, a moduli space of vector
bundles over $\mathbb{P}^{2}$ is also a moduli space of
representations for this quiver with relations. In fact, one can
show that the moduli spaces of representations that occur in this way
may be taken to parametrise certain rather special representations
and using the known results on moduli spaces of representations of
quivers it is possible to show that these moduli spaces are birational
to a suitable number of matrices up to simultaneous conjugacy.

The next section introduces the notation that we shall use in this
paper. Then section \ref{leftgeneral} studies particular moduli spaces
of representations that include all the examples we shall need. In
section \ref{modp2} we show how to reduce a moduli space of vector
bundles over $\mathbb{P}^{2}$ to one of these moduli spaces of
representations.

\section{Terminology}
\label{term}

We introduce some notation and terminology. The terminology we
introduce for representations of a multiplication below is useful
since this is a particular case of representations of a quiver with
relations which can reasonably be presented in a basis-free way.  Let
$U$ be a finite dimensional vector space. A \definition{representation
$R$ of the vector space $U$} is a triple $(R(0),R(1),R(\phi))$ where
$R(0)$ and $R(1)$ are finite dimensional vector spaces and
$R(\phi)\colon R(0)\otimes U\rightarrow R(1)$ is a linear map. Its
\definition{dimension vector} is $\ddim R=(\dim R(0)\ \dim R(1))$. The
representations of dimension vector $\alpha=(a\ b)$ are parametrised
by the vector space $R(U,\alpha)=\Hom(k^{a}\otimes V,k^{b})$ on
which the algebraic group $Gl_{a}\times Gl_{b}$ acts by change of
bases. The orbits correspond to the isomorphism classes of
representations. Of course the category of representations of a vector
space of dimension $n$ is just the category of representations of a
quiver with two vertices and $n$ arrows from the first to the second
vertex as one sees by choosing a basis of $U$. We note that
$R(\phi)\dual\colon R(1)\dual\rightarrow R(0)\dual\otimes U\dual$
gives a representation of $U$, $R\dual$, from the linear map
$R\dual(\phi)\colon R(1)\dual\otimes U\rightarrow R(0)\dual$.

A \definition{multiplication} is a quadruple $(U,V,W,f)$ where $U$,
$V$ and $W$ are vector spaces and $f\colon U\otimes V\rightarrow W$ is
a linear map.  We shall usually talk of the multiplication $f$. A
\definition{representation $R$ of the multiplication} $f$ is a
sextuple $(R(0),R(1),R(2),R(\phi_{01}),R(\phi_{12}),R(\phi_{02}))$
where each $R(i)$ is a finite dimensional vector space and
$R(\phi_{01})\colon R(0)\otimes U\rightarrow R(1)$,
$R(\phi_{12})\colon R(1)\otimes V\rightarrow R(2)$ and
$R(\phi_{02})\colon R(0)\otimes W\rightarrow R(2)$ are linear maps such that
\begin{equation}
\label{e0}
(\phi_{01}\otimes I_{V})\phi_{12}= (I_{R(0)}\otimes f)\phi_{02} 
\end{equation}
as linear maps from $R(0)\otimes U\otimes V$ to $R(2)$. We shall
eventually be interested in representations of the multiplication
$\sigma\colon k^{3}\otimes k^{3}\rightarrow S^{2}(k^{3})$ since the
results of Beilinson \cite{Beilinson} show that the derived category
of coherent sheaves on $\mathbb{P}^{2}$ is equivalent to the derived
category of representations of the multiplication $\sigma$.  The
\definition{dimension vector} of a representation $R$ of the
multiplication $f$ is $\ddim R=(\dim R(0)\ \dim R(1)\ \dim R(2))$. The
representations of dimension vector $\alpha=(a\ b\ c)$ are
parametrised by the closed subvariety $R(f,\alpha)$ of
$\Hom(k^{a}\otimes U,k^{b})\times \Hom(k^{b}\otimes
V,k^{c})\times \Hom(k^{a}\otimes W,k^{c})$ of triples
$(\phi_{01},\phi_{12},\phi_{02})$ satisfying equation \ref{e0}. The
algebraic group $Gl_{\alpha}=Gl_{a}\times Gl_{b}\times Gl_{c}$ acts
via change of bases on $R(f,\alpha)$ and the orbits correspond to the
isomorphism classes of representations of dimension vector $\alpha$ of
the multiplication. Vector space duality gives rise to linear maps
$R(2)\dual\otimes V\rightarrow R(1)\dual$, $R(1)\dual\otimes
U\rightarrow R(0)\dual$ and $R(2)\dual\otimes W\rightarrow R(0)\dual$
and this gives a representation of the multiplication
$\hat{f}\colon V\otimes U\rightarrow W$ obtained from $f$ by switching $U$
and $V$.

$R(f,\alpha)$ is in general a reducible variety and the description of
its components and their orbit spaces in complete generality is not
something we shall undertake in this paper. However, in the case
arising from vector bundles over $\mathbb{P}^{2}$, we can restrict to
components parametrising representations of a fairly nice form for
which the moduli space is relatively comprehensible. 

Given a dimension vector $\alpha$ for a multiplication, $\alpha_{l}$
will be the dimension vector $(\alpha(0)\ \alpha(1))$ and $\alpha_{r}$
will be the dimension vector $(\alpha(1)\ \alpha(2))$. If $R$ is a
representation of the multiplication $f\colon U\otimes V\rightarrow
W$, then $\alpha_{l}$ is the dimension vector of the representation of
the vector space $U$ obtained by restriction and $\alpha_{r}$ is the
dimension vector of the representation of the vector space $V$
obtained by restriction. Thus we have morphisms from $R(f,\alpha)$ to
$R(U,\alpha_{l})$ and to $R(V,\alpha_{r})$. Let $C$ be an irreducible
component of $R(f,\alpha)$ such that the morphism to $R(U,\alpha)$ is
dominant; then we shall say that $C$ is a \definition{left general}
component and that a general representation in the component $C$ is
\definition{left general}. The term \definition{right general} is
defined in a similar way using the morphism to $R(V,\alpha_{r})$.
Later results about left general representations have analogous
results for right general representations since duality will carry the
one to the other. 

Representations of a vector space or of a multiplication are special
cases of the more general notion of representations of a quiver or of
a quiver with relations as one sees by choosing bases for the vector
spaces $U$, $V$ and $W$. We refer the reader to section 2 of
\cite{qmod} for the terminology we shall use for representations of a
quiver.

A relation on a quiver $Q$ is a linear combination of paths
$r=\sum_{i=1}^{n} \lambda_{i}p_{i}$ such that the initial and terminal
vertex of the path $p_{i}$ are all equal. Given a representation $R$
of the quiver,we may extend our notation by defining
$R(p)=R(a_{1})\dots R(a_{n})$ if $p=a_{1}\dots a_{n}$ and
$R(r)=\sum_{i=1}^{n}\lambda_{i}R(p_{i})$.  Given a set of relations
$I$, a representation $R$ of the quiver with relations $(Q,I)$ is a
representation of the quiver $Q$ such that $R(r)=0$ for every $r\in
I$. The category $\Rep(Q,I)$ is the full subcategory of $\Rep(Q)$
whose objects are the representations of the quiver with relations
$(Q,I)$.

A family $\mathcal{R}$ of representations of dimension vector $\alpha$
of the quiver with relations $(Q,I)$ over an algebraic
variety $X$ is a collection of vector bundles $R(v)$ for each vertex
$v$ and homomorphisms of vector bundles $\mathcal{R}(a)\colon
\mathcal{R}(ia)\rightarrow \mathcal{R}(ta)$ for each arrow $a$ such
that for each relation $r$, $\mathcal{R}(r)=0$. Given a point $p\in X$
there is an associated representation $\mathcal{R}_{p}$, the fibre of
$\mathcal{R}$ above the point $p$. When $X$ is an irreducible
algebraic variety we shall say that the family is
\definition{irreducible} and in this case we shall say that a
representation $R$ \definition{is of type} $\mathcal{R}$ if there
exists a point $p$ such that $R\cong \mathcal{R}_{p}$. Again when
$\mathcal{R}$ is an irreducible family, we shall say that a general
representation of type $\mathcal{R}$ has property $P$ if there exists
a dense open subvariety $O$ of $X$ such that for all $p\in O$,
$\mathcal{R}_{p}$ has property $p$. 

Given a dimension vector $\alpha$, the vector space $R(Q,\alpha)=
\oplus_{a\in A} {^{\alpha(ia)}k^{\alpha(ta)}}$ parametrises the
representations of the quiver of dimension vector $\alpha$ and this
carries a family $\mathcal{R}$ of representations of the quiver which
we shall refer to as the \definition{canonical} family. If $V$ is a
locally closed subvariety of $R(Q,\alpha)$, then the restriction of
the canonical family to $V$ will be called the \definition{canonical
family on $V$}.  There is a closed subvariety $R(Q,I,\alpha)$, the
\definition{representation space} of dimension vector $\alpha$ for the
quiver with relations $(Q,I)$ consisting of those points $p$ such that
$\mathcal{R}_{p}(r)=0$ for every relation $r$ in $I$. In general this
is a reducible algebraic variety. 

The algebraic group $Gl_{\alpha}= \times_{v}Gl_{\alpha(v)}$ acts on
this family compatibly with an action of
$PGl_{\alpha}=Gl_{\alpha}/k^{*}$ on $R(Q,I,\alpha)$ and the orbits of
$PGl_{\alpha}$ on $R(Q,I,\alpha)$ correspond to the isomorphism
classes of representations of dimension vector $\alpha$. In addition,
if $p\in R(Q,\alpha)$, the stabiliser of $p$ in $Gl_{\alpha}$ acts on
the representation defined by the point $p$ as the units of its
endomorphism ring. 

Let $X$ be an algebraic variety on which the algebraic group $G$
acts. Let
\begin{equation*}
1\rightarrow k^{*}\rightarrow \tilde{G}\rightarrow G\rightarrow 1
\end{equation*}
be a short exact sequence of algebraic groups.  Let $E$ be a vector
bundle over $X$ on which $\tilde{G}$ acts compatibly with the action
of $G$. Then $k^{*}$ acts on the fibres of $E$ and if this action is
via the character $\phi_{w}(\lambda)=\lambda^{w}$ then we shall say
that $E$ is \definition{a $\tilde{G}$ vector bundle of weight} $w$. A
morphism of $\tilde{G}$ vector bundles of weight $w$ is a morphism of
vector bundles that is also $\tilde{G}$ equivariant. 
We see that if $\mathcal{R}$ is a canonical family of representations
of dimension vector $\alpha$ then each $\mathcal{R}(v)$ is a
$Gl_{\alpha}$ vector bundle of weight $1$.

Two irreducible families of representations of the quiver with
relations $(Q,I)$, $\mathcal{R}$ over $X$ and $\mathcal{S}$ over $Y$,
are said to be \definition{birationally representation equivalent} if
there exist open dense subvarieties $O\subset X$ and $O'\subset Y$ such
that for all $p\in O$, there exists $q\in Y$ such that
$\mathcal{R}_{p}\cong \mathcal{S}_{q}$ and for all $q\in O'$ there
exists $p\in X$ such that $\mathcal{R}_{p}\cong \mathcal{S}_{q}$. We
shall say that the family $\mathcal{S}$ \definition{birationally
contains} the family $\mathcal{R}$ if there exists an open subvariety
$O\subset X$ such that for all $p\in O$, there exists $q\in Y$ such
that $\mathcal{R}_{p}\cong \mathcal{S}_{q}$. We shall say that the
family $\mathcal{R}$ is \definition{birationally constant} if there
exists a dense open subvariety $O$ of $X$ such that
$\mathcal{R}_{p}\cong \mathcal{R}_{q}$ for all $p$ and $q$ in $O$.

We shall say that a family $\mathcal{R}$ of representations of
dimension vector $\alpha$ is \definition{reducible to matrix normal
form of type $h$} if it is birationally representation equivalent to a
family $\mathcal{S}$ over an irreducible algebraic variety $X$ on
which $PGl_{h}$ acts such that $X$ is $PGl_{h}$ birational to
$M_{h}(k)^{t}$ for some non-negative integer $t$ where $PGl_{h}$ acts
by conjugation on each factor, and each $\mathcal{S}(v)$ is a $Gl_{h}$
vector bundle of weight $1$ so that the action of the stabiliser in
$Gl_{h}$ of a point $p$ in $X$ acts on $\mathcal{S}_{p}$ as the unit
group of the endomorphism ring of $\mathcal{S}_{p}$. It is said to be
\definition{reducible to matrix normal form} if in addition
$h=\hcf_{v}(\alpha_{v})$. A family $\mathcal{S}$ with these properties
will be called an $\mathcal{R}$-\definition{standard} family.

\section{Left general components}
\label{leftgeneral}

The purpose of this section is to show that a moduli space of left
general representations of a multiplication is birational to a
suitable number of matrices up to simultaneous conjugacy when there is
a left general representation with trivial endomorphism ring. We shall
actually prove a more general result about families of representations
of a quiver with relations which may be of future use. We begin
by showing that there is only one left general component of given
dimension vector. 

\begin{lemma}
\label{uniqlg}
Let $\alpha$ be a dimension vector for the multiplication $f$.
There is a unique left general component of dimension vector $\alpha$ .
\end{lemma}
\begin{proof}
Let $\alpha=(a\ b\ c)$ and let $f\colon U\otimes V\rightarrow W$ be
the multiplication. Let $O$ be the open subvariety of $R(U,(a\ b))$ of
points such that the linear map $I\otimes f\oplus -R_{p}(\phi)\otimes
I\colon k^{a}\otimes U\otimes V\rightarrow k^{a}\otimes W\oplus
k^{b}\otimes V$ has maximal rank. Then $O$ carries a vector bundle $E$
whose fibre above the point $p$ is simply the cokernel of $I\otimes
f\oplus -R_{p}(\phi)\otimes I$. We consider the vector bundle
$E\dual\otimes k^{c}$ which is an irreducible variety and has a
morphism to $R(f,(a\ b\ c))$ whose image must contain an open dense
subvariety of any left general component since its image contains
every point $q$ such that the left restriction of $R_{q}$ is
isomorphic to a representation of the form $R_{p}$ for $p\in O$. Its
image is irreducible and consequently must lie in one of these
components so it follows that there must be precisely one left general
component and the image of $E\dual\otimes k^{c}$ must lie in this
unique component.
\end{proof}

We shall need to be able to recognise that a representation of a
multiplication lies in the left general component for its dimension
vector. 

\begin{lemma}
\label{reclg}
Let $f\colon U\otimes V\rightarrow W$ be a multiplication such that
$f$ is surjective. Let $K$ be the kernel of $f$. Let $R$ be a
representations of the multiplication $f$ such that the composition 
of the linear maps from $R(0)\otimes K$ to $R(0)\otimes U\otimes V$ to
$R(1)\otimes V$ is injective. Then $R$ lies in the left general
component of representations of its dimension vector.
\end{lemma}
\begin{proof}
The injectivity of this linear map implies that the map considered in
the proof of the previous lemma has maximal rank.
\end{proof}

In order to prove the main theorem of this section we shall need to
summarise the information we already have on representations of a
vector space. 

We consider representations of the vector space $U$ of dimension
vector $(a\ b)$. Let $g=\hcf(a,b)$. After choosing a basis of $U$, we
see that we are simply considering representations of a generalised
Kronecker quiver, that is, a quiver $Q(u)$ where $u=\dim U$ which has
two vertices $v$ and $w$ and $u$ arrows from $v$ to $w$. We shall be
interested in the structure of general representations of these
quivers. If $u=1$, nothing difficult happens.  For a dimension vector
$(a\ b)$, if $a<b$, the general representation is isomorphic to $(k\
k)^{a}\oplus (0\ k)^{b-a}$; if $a=b$, the general representation is
isomorphic to $(k\ k)^{a}$ and if $a>b$ the general representation is
isomorphic to $(k\ 0)^{a-b}\oplus (k\ k)^{b}$. For $u>1$, we need to
introduce some terminology. The results in this case which are
described below may be found on page 159 in \cite{Kac}. The projective
indecomposable representations of $U$ are $P_{0}=(0\ k)$ and
$P_{1}=(k\ U)$ whilst the injective indecomposable representations are
$I_{0}=(k\ 0)$ and $I_{1}=(U\dual\ k)$ which are simply the dual
representations to the projective representations. We assume that we
have constructed representations $\{P_{i}:i=0\rightarrow t\}$ and we
have shown that $\Hom(P_{i},P_{i+1})\cong U$ if $i$ is even whilst
$\Hom(P_{i},P_{i+1})\cong U\dual$ if $i$ is odd. Then if $t$ is
odd, there is a canonical homomorphism from $P_{t-1}$ to $U\dual
\otimes P(t)$ and we define $P_{t+1}$ to be the cokernel whilst if $t$
is even, there is a canonical homomorphism from $P_{t-1}$ to $U\otimes
P_{t}$ whose cokernel we define to be $P_{t+1}$. One may check that
the inductive hypothesis has been extended and thus we have defined
representations $P_{n}$ for all integers $n$. These representations
are called the preprojective representations of $U$. We define
$I_{n}={P_{n}}\dual$ and these are the preinjective
representations. In the case where $u=2$, $P_{i}\cong (S^{i-1}(U)\
S^{i}(U))$ for $i>0$ where the linear map from $S^{i-1}(U)\otimes U$
to $S^{i}(U)$ is the obvious one. These representations are important
to us because they allow us to describe the general representations of
arbitrary dimension vector for a vector space $U$.

We first state the results for $u=2$. We consider the dimension vector
$(a\ b)$. If $a\neq b$, then $Gl_{a}\times Gl_{b}$ has an open orbit
in $R(U,(a\ b))$ and the corresponding representation will be called
$G(a\ b)$. If $a<b$ then for some integer $m$, $\tfrac{m-1}{m}\leq
\tfrac{a}{b}<\tfrac{m}{m+1}$. If $\tfrac{a}{b}=\tfrac{m-1}{m}$ then
$G(a\ b)$ is isomorphic to ${P_{m-1}}^{b/m}$ and clearly
$\tfrac{b}{m}=\hcf(a,b)$. Otherwise $G(a\ b)$ is isomorphic to
$P_{m-1}^{c}\oplus P_{m}^{d}$ where $c$ and $d$ are non-zero integers
such that $\hcf(c,d)=\hcf(a,b)$. If $a>b$ then duality leads to the
same picture using preinjective representations. Thus if
$\tfrac{b}{a}=\tfrac{m-1}{m}$ then $G(a\ b)\cong {I_{m-1}}^{a/m}$ where
$\tfrac{a}{m}=\hcf(a,b)$ and if
$\tfrac{m}{m+1}>\tfrac{b}{a}>\tfrac{m-1}{m}$ then $G(a\ b)$ is
isomorphic to $I_{m-1}^{c}\oplus I_{m}^{d}$ where
$\hcf(c,d)=\hcf(a,b)$. If $a=b$, then $Gl_{a}\times Gl_{b}$ does not
have an open orbit on $R(U,(a\ b))$, however, this case is essentially
one $a$ by $a$ matrix up to simultaneous conjugacy as we shall now
see. There is a family of representations of the vector space $U$ on
the algebraic variety $X=M_{a}(k)$ on which $PGl_{a}$ acts by
conjugation defined as follows. Let $\{u_{1},u_{2}\}$ be a basis of
$U$; let $\mathcal{R}(0)=\mathcal{R}(1)=k^{a}\times X$ and
$\mathcal{R}_{p}(\phi)(v_{1}\otimes u_{1}+v_{2}\otimes u_{2})=
v_{1}+v_{2}p$ (recall that $p\in X=M_{a}(k)$). Note that this is a
reduction to matrix normal form for the dimension vector $(a\ a)$.

This leaves the case where $u>2$. If $(a\ b)$ is a Schur root then
$(a\ b)$ is reducible to matrix normal form. Otherwise $Gl_{a}\times
Gl_{b}$ has an open orbit on $R(U,(a\ b))$ and the corresponding
representation $G(a\ b)$ is isomorphic to $P_{m}^{c}\oplus
P_{m+1}^{d}$ for suitable integers $m,c\text{ and }d$ when $a<b$ and
is isomorphic to $I_{m}^{c}\oplus I_{m+1}^{d}$ when $a>b$. Further
$\hcf(c,d)=\hcf(a,b)$ which is not demonstrated in these theorems but
may be checked quickly by use of the reflection functors to reduce to
the case where $m=0$ where it is clear.

In all the cases where $G(a\ b)$ is defined and is isomorphic to
either $P_{m}^{g}$ or $I_{m}^{g}$ the dimension vector is reducible to
matrix normal form trivially; so $X$ is a point on
which $PGl_{g}$ acts trivially; the family is simply the
representation $P_{m}^{g}$ which has endomorphism ring $M_{g}(k)$ and
$Gl_{g}$ acts on the family via its action as the group of
automorphisms of $P_{m}^{g}$.

Thus we have a useful dichotomy; either the dimension vector $(a\ b)$
is reducible to matrix normal form or else $Gl_{a}\times Gl_{b}$ has
an open orbit on $R(U,(a\ b))$ and the corresponding representation
$G(a\ b)$ is isomorphic to either $P_{m}^{c}\oplus P_{m+1}^{d}$ for
positive integers $c$ and $d$ and non-negative integer $m$ or else to
$I_{m}^{c}\oplus I_{m+1}^{d}$ with the same conditions on $c$, $d$ and
$m$. In these cases, $\hcf(c,d)=\hcf(a,b)$. Further, $P_{m}^{c}\oplus
P_{m+1}^{d}$ has a unique subrepresentation of dimension vector $\ddim
P_{m+1}^{d}$ whilst $I_{m}^{c}\oplus I_{m+1}^{d}$ has a unique
subrepresentation of dimension vector $\ddim I_{m}^{c}$. These
statements remain true for $u=1$ by setting $P_{0}=(0\ k)$, $P_{1}=(k\
k)=I_{1}$ and $I_{0}=(k\ 0)$. We note that the endomorphism ring of
$P_{m}\oplus P_{m+1}$ or of $I_{m}\oplus I_{m+1}$ is isomorphic to
$\left(\begin{smallmatrix} k &U''\\ 0 &k
\end{smallmatrix}\right)
$ where $U''$ is isomorphic to $U$ or $U\dual$. The corresponding
statement remains true when $u=1$. Therefore, if $(a\ b)$ is not
reducible to matrix normal form, the representation $G(a\ b)$ is
defined and its endomorphism ring is Morita equivalent to 
$\left(\begin{smallmatrix} k &U\\ 0 &k
\end{smallmatrix}\right)
$. It will be useful to have a common notation for these two subcases
of the second case. We shall say that $G(a\ b)\cong G_{1}^{a'}\oplus
G_{2}^{b'}$ where $a'$ and $b'$ are positive integers such that
$\hcf(a',b')=\hcf(a,b)$ and $G(a\ b)$ has a unique subrepresentation of
dimension vector $a'\ddim G_{1}=(d\ e)$ where $\hcf(d,e)=a'$.  

For the next result we shall regard $Gl_{h}$ as a subgroup of
$Gl_{\alpha}$ via the diagonal embedding in each factor
$Gl_{\alpha(v)}$ and $PGl_{h}$ as the corresponding subgroup of
$PGl_{\alpha}$.

\begin{lemma}
\label{standincan}
Let $C$ be an irreducible component of $R(Q,I,\alpha)$ and let
$\mathcal{R}$ be the canonical family of representations of dimension
vector $\alpha$ on $C$. Assume that $\mathcal{R}$ is reducible to
matrix normal form. Then $C$ has a $PGl_{h}$ equivariant subvariety
$Y$ such that the restriction of $\mathcal{R}$ to $Y$ is a standard
$\mathcal{R}$-family and therefore $C$ is $PGl_{h}$ birational to
$M_{h}(k)^{s}\times^{PGl_{h}}PGl_{\alpha}$ for some integer $s$. In
particular, a moduli space of representations of type $C$ is
birational to a suitable number of matrices up to simultaneous
conjugacy.
\end{lemma}
\begin{proof}
Let $\mathcal{S}$ be a $\mathcal{R}$-standard family over the
algebraic variety $X$. Then for each vertex $v$, $\mathcal{S}(v)$ is a
vector bundle of weight $1$ for $Gl_{h}$ of rank $\alpha(v)=h\beta(v)$
for a suitable dimension vector $\beta$. Then by the local isomorphism
theorem, lemma 3.1 of \cite{qmod}, there exists an open $PGl_{h}$
equivariant subvariety $X'$ of $X$ such that the restriction of
$\mathcal{S}(v)$ to $X'$ is isomorphic as $Gl_{h}$ vector bundle to
$k^{\beta(v)}\otimes k^{h}\times X'$ where $Gl_{h}$ acts trivially on
$k^{\beta(v)}$ and diagonally on $k^{h}\times X'$. We choose bases for
each $k^{\beta(v)}$ and a basis of $k^{h}$ which give bases for each
$k^{\beta(v)}\otimes k^{h}\cong k^{\alpha(v)}$ and determines a
homomorphism from $Gl_{h}$ to each $Gl_{\alpha(v)}$ and hence to
$Gl_{\alpha}$ which we may regard as a diagonal embedding. Each arrow
$a$ determines a morphism of vector bundles from $k^{\beta(ia)}\otimes
k^{h}\times X'$ to $k^{\beta(ta)}\otimes k^{h}\times X'$ and hence
determines a morphism of algebraic varieties from $X'$ to
${^{\alpha(ia)}k^{\alpha(ta)}}$ and hence we have a morphism of
algebraic varieties from $X'$ to $R(Q,\alpha)$ which is $PGl_{h}$
equivariant and has the property that the pullback of the canonical
family on the image of the morphism is $\mathcal{S}$. It follows that
the image actually lies in $R(Q,I,\alpha)$ and in fact must lie in the
component $C$. Since the stabiliser in $Gl_{h}$ of a point $p$ in $X'$
is isomorphic to the units of $\End(\mathcal{S}_{p})$ and so is the
stabiliser in $Gl_{\alpha}$ of its image $q$ in $C$, we deduce that
the morphism from $X'$ to $C$ is injective; indeed that the morphism
from $X'\times^{PGl_{h}}PGl_{\alpha}$ to $C$ is injective and this
latter morphism is also dominant since $\mathcal{S}$ is birationally
representation equivalent to $\mathcal{R}$.
\end{proof}

Given irreducible families $\mathcal{R}$ over $X$ and $\mathcal{S}$
over $Y$, the two functions on $X\times Y$ that assign to the point
$(p,q)$ the values $\hom(\mathcal{R}_{p},\mathcal{S}_{q})$ and
$\ext(\mathcal{R}_{p},\mathcal{S}_{q})$ are upper semicontinuous and
consequently there exists an open dense subvariety $O\subset X\times
Y$ where these functions are constant and minimal; these minimal
values we shall call $\hom(\mathcal{R},\mathcal{S})$ and
$\ext(\mathcal{R},\mathcal{S})$. If $X$ is a point and $R$ is the
corresponding representation, we write $\hom(R,\mathcal{S})$ and
$\ext(R,\mathcal{S})$; $\hom(\mathcal{R},S)$ and $\ext(\mathcal{R},S)$
are defined similarly for a representation $S$.

Let $O$ be the dense open subvariety of $X\times Y$ consisting of
points $(p,q)$ where
$\ext(\mathcal{R}_{p},\mathcal{S}_{q})=\ext(\mathcal{R},\mathcal{S})$.
Then there is a vector bundle $E(\mathcal{R},\mathcal{S})$ over $O$
whose fibre above the point $(p,q)$ is
$\Ext(\mathcal{R}_{p},\mathcal{S}_{q})$ and there is a family of
representations over $E(\mathcal{R},\mathcal{S})$ of extensions of
representations of type $\mathcal{R}$ on representations of type
$\mathcal{S}$. We shall call this family \definition{the extension
family of $\mathcal{R}$ on $\mathcal{S}$},
$\mathcal{E}(\mathcal{R},\mathcal{S})$. If $Y$ is just a point then
$\mathcal{S}$ is just a representation $S$ and we shall refer to the
extension family of $\mathcal{R}$ on $S$,
$\mathcal{E}(\mathcal{R},S)$; similarly we define the extension family
of $R$ on $\mathcal{S}$ for a representation $R$.

Let $S$ be a representation such that $\End(S)=k$ and
$\Ext(S,S)=0$. We shall say that a representation $R$ has an $S$-socle
if the natural map from $\Hom(S,R)\otimes S$ to $R$ is injective and
the image will be called the \definition{$S$-socle}. If $\Hom(S,R)=0$
we shall say that $R$ is \definition{left $S$-free}. When $R$ has an
$S$-socle, $T$, then it is clear that $R/T$ is left $S$-free (apply
$\Hom(S,\ )$ to the short exact sequence $\ses{T}{R}{R/T}$). Dually we
say that $R$ has an $S$-top if the natural map from $R$ to
$\Hom(R,S)\dual\otimes S$ is surjective; if $K$ is the kernel of this
homomorphism we say that $R/K$ is the \definition{$S$-top} of $R$;
again it follows that $\Hom(K,S)=0$ and we define a representation $R$
to be \definition{right $S$-free} if $\Hom(R,S)=0$.

Now let $\mathcal{R}$ be an irreducible family over the algebraic
variety $X$. Let $O$ be the dense open subvariety on which
$\hom(S,\mathcal{R}_{p})=c=\hom(S,\mathcal{R})$. Assume that there is
a point $p$ in $O$ such that $\mathcal{R}_{p}$ has an $S$-socle and so
the natural map from $\Hom(S,\mathcal{R}_{p})\otimes S$ to
$\mathcal{R}_{p}$ is injective. Then a general representation of type
$\mathcal{R}$ must have an $S$-socle which is isomorphic to
$S^{c}$. Then on a suitable dense open subvariety $U$ of $X$, we have
the \definition{associated left $S$-free family}, $\mathcal{R}'$,
where the representation $\mathcal{R}'_{p}$ is $\mathcal{R}_{p}/T_{p}$
where $T_{p}\cong S^{c}$ is the $S$-socle of $\mathcal{R}_{p}$. The
family $\mathcal{R}$ is birationally contained in
$\mathcal{E}(\mathcal{R}',S^{c})$ and we shall say that $\mathcal{R}$
is \definition{left general with respect to} $S$ if they are
birationally representation equivalent. We note that if $\mathcal{R}$
is the canonical family on a component $C$ of $R(Q,I,\alpha)$ it is
forced to be left general with respect to $S$ whenever a general
representation of type $C$ has an $S$-socle.

Dually, a general representation of type $\mathcal{R}$ may have an
$S$-top isomorphic to $S^{d}$; then we have a family on $X$,
$\mathcal{R}''$, which we call the associated right $S$-free family,
where $\mathcal{R}''_{p}$ is the kernel of the surjection from
$\mathcal{R}_{p}$ onto $S^{d}$ and we define $\mathcal{R}$ to be
\definition{right general with respect to} $S$ if $\mathcal{R}$ and
$\mathcal{E}(S^{d},\mathcal{R}'')$ are birationally representation
equivalent.

The following lemma is our main reduction.

\begin{lemma}
\label{mnfontop}
Let $\mathcal{R}$ be an irreducible family of representations of
dimension vector $\alpha$ of the quiver with relations $(Q,I)$ over
the algebraic variety $X$ such that a general representation of type
$\mathcal{R}$ has trivial endomorphism ring. Let $S$ be a
representation such that $\End(S)=k$ and $\Ext(S,S)=0$. Assume that a
general representation of type $\mathcal{R}$ has an $S$-socle and that
$\mathcal{R}$ is left general with respect to $S$. Let
$c=\hom(S,\mathcal{R})$. Assume that the associated left $S$-free
family is reducible to matrix normal form of type $g$.  Then
$\mathcal{R}$ is reducible to matrix normal form of type
$\hcf(g,c)$. Similarly, if a general representation of type
$\mathcal{R}$ has an $S$-top, $\mathcal{R}$ is right general with
respect to $S$, the associated right $S$-free family is reducible to
matrix normal form of type $h$ and $\hom(\mathcal{R},S)=d$ then
$\mathcal{R}$ is reducible to matrix normal form of type $\hcf(h,d)$.
\end{lemma}
\begin{proof}
This argument essentially occurs as a special case of the proof of
theorem 6.1 in \cite{qmod} to which the reader should refer for greater
detail.

We deal only with the first case since the second case has the same
proof. Let $\mathcal{R}'$ be the associated left $S$-free family on
$X$ (after replacing $X$ by an open dense subvariety). Let
$\mathcal{S}$ be an $\mathcal{R}'$-standard family over $Y$ where $Y$
is $PGl_{g}$ equivariant to a dense open subvariety of $M_{g}(k)^{s}$
for some integer $s\geq 0$. Let $E=E(\mathcal{S},S^{c})$ and let
$\mathcal{E}=\mathcal{E}(\mathcal{S},S^{c})$ be the extension family
of $\mathcal{S}$ on $S^{c}$. Then since we assume that $\mathcal{R}$
is left general with respect to $S$, $\mathcal{R}$ and $\mathcal{E}$
are birationally representation equivalent. 
 
Let $\ext(\mathcal{S},S)=t$. After shrinking $Y$ a little we may
assume that for all $p\in Y$, $\ext(\mathcal{S}_{p},S)=t$. Let
$\beta=(g\ c)$ be a dimension vector for the quiver $Q'$ which has two
vertices $v$ and $w$, $s$ arrows from $v$ to itself and $t$ arrows
from $v$ to $w$. Then $PGl_{\beta}$ acts on $E$ whilst $Gl_{\beta}$
acts on $\mathcal{E}$ so that $\mathcal{E}(v)$ is a $Gl_{\beta}$
vector bundle of weight $1$ and $E$ is $PGl_{\beta}$ birational to
$R(Q',\beta)$. 

Since $\mathcal{R}$ and $\mathcal{E}$ are birationally representation
equivalent, a general representation of type $\mathcal{E}$ has an
$S$-socle which must coincide with its obvious subrepresentation
isomorphic to $S^{c}$. We therefore pass to a dense open $PGl_{\beta}$
subvariety of $E$ where this is true. Then the orbits of $PGl_{\beta}$
correspond to the isomorphism classes of representations in the
restriction of $\mathcal{E}$ to this subvariety. Since a general
representation of type $\mathcal{R}$ and hence of type $\mathcal{E}$
has trivial endomorphism ring, it follows that $PGl_{\beta}$ has
trivial stabilisers generically on $E$ and hence $\beta$ is a Schur
root for the quiver $Q'$. Thus the main result of \cite{qmod}, theorem
6.3, allows us to conclude that $\beta$ is reducible to matrix normal
form and hence by lemma \ref{standincan} there exists a $PGl_{m}$
equivariant subvariety $Z$ of $R(Q',\beta)$ where $m=\hcf(g,c)$ such
that the restriction of the canonical family on $R(Q',\beta)$ to $Z$
is standard and $ZPGl_{\beta}$ is a dense subvariety of
$R(Q',\beta)$. Since $E$ and $R(Q',\beta)$ are $PGl_{\beta}$
birational, there is a corresponding $PGl_{m}$ equivariant subvariety
$Z'$ of $E$ and the restriction of $\mathcal{E}$ to $Z'$ is what we
want.
\end{proof}

Now let $S$ and $S_{1}$ be representations such that
$\End(S_{1})=k=\End(S)$, $\Ext(S_{1},S_{1})=0=\Ext(S,S)$ and also
$\Hom(S,S_{1})=0=\Hom(S_{1},S)=\Ext(S,S_{1})$. Let $\E(S,S_{1})$ be
the full subcategory of representations that contains $S$ and $S_{1}$
and is closed under extensions; we call the pair $(S,S_{1})$ a
\definition{Kronecker reduction pair}. Then following Ringel
\cite{Ringel} we have the following lemma.

\begin{lemma}
\label{ringelred}
Let $S$ and $S_{1}$ be a Kronecker reduction pair. Then every
representation in $\E(S,S_{1})$, the full subcategory of
representations that contain $S$ and $S_{1}$ and is closed under
extensions, has an $S$-socle such that the factor is isomorphic to
$S_{1}^{n}$ for some integer $n$. Let $t=\ext(S_{1},S)$ and let 
\begin{equation*}
\ses{S^{t}}{S'}{S_{1}} 
\end{equation*}
be the canonical extension of $S_{1}$ on $S^{t}$. Then $\Ext(S\oplus
S',M)=0$ for all $M$ in $\E(S,S_{1})$, and $\Hom(S\oplus S',\ )$
induces a natural equivalence with the category of representations of
the $t$th Kronecker quiver whose inverse is given by $\ \otimes
(S\oplus S')$. This functor also induces isomorphisms on Ext groups.
\end{lemma}
\begin{proof}
If $M$ and $N$ are representations for which we have short exact
sequences 
\begin{align*}
&\ses{S^{m_{1}}}{M}{S_{1}^{n_{1}}}\\
&\ses{S^{m_{2}}}{N}{S_{1}^{n_{2}}} 
\end{align*}
and a short exact sequence
\begin{equation*}
\ses{M}{L}{N} 
\end{equation*}
then it is clear that the induced extension of the subrepresentation
$S^{m_{2}}$ over $M$ splits; thus $L$ has a subrepresentation $T$ 
isomorphic to $S^{m_{1}+m_{2}}$ such that the factor is isomorphic to
$S_{1}^{n_{1}+n_{2}}$ and $T$ must be the $S$-socle of $L$. Thus by
induction all objects in $\E(S,S_{1})$ have the required structure.

By construction, $\Ext(S',S)=0$ and clearly $\Ext(S',S_{1})=0$;
therefore $\Ext(S\oplus S',M)=0$ for all $M$ in $\E(S,S_{1})$ as
required. Now let $M$ be any object in $\E(S,S_{1})$; so there is a
short exact sequence 
\begin{equation*}
\ses{S^{m}}{M}{S_{1}^{n}}.
\end{equation*}
Then $\hom(S,M)=m$ and since $\ext(S',S)=0$, $\hom(S',M)=n$ so that
$\Hom(S\oplus S',M)$ is a representation of dimension vector $(n\ m)$
of the $t$th Kronecker quiver. In fact, the natural homomorphism from
$\Hom(S,M)\otimes S\oplus \Hom(S',M)\otimes S$ to $M$ is surjective
with kernel isomorphic to $S^{s}$ for some integer $s$. Using this
short exact sequence it is a simple matter to check that the two
functors are mutually inverse and that $\Ext(M,N)$ is preserved by
this functor.
\end{proof}

Let $\mathcal{R}$ be some irreducible family of representations of the
quiver with relations $(Q,I)$ over the algebraic variety $X$. We shall
say that $\mathcal{R}$ has a \definition{Kronecker reduction of type}
$(S,S_{1})$ to the dimension vector $(a\ b)$ for the $t$th Kronecker
quiver if there exists a Kronecker reduction pair $(S,S_{1})$ where
$\ext(S_{1},S)=t$ such that a general representation of type
$\mathcal{R}$, $R$, has a subrepresentation isomorphic to $S^{b}$ with
factor isomorphic to $S_{1}^{a}$ for suitable integers $b$ and $a$ (so
a general representation of type $\mathcal{R}$ has an $S$-socle) and
$\mathcal{R}$ is left general with respect to $S$. Note that this is a
self-dual condition since being left general with respect to $S$ is
equivalent to being right general with respect to $S_{1}$.

We note the following consequence of lemma \ref{ringelred}.

\begin{lemma}
\label{kr}
Let $\mathcal{R}$ be an irreducible family of representations of
dimension vector $\alpha$ of the quiver with relations $(Q,I)$ over
the algebraic variety $X$. Assume that $\mathcal{R}$ has a Kronecker
reduction of type $(S,S_{1})$ to the dimension vector $(a\ b)$ for the
$t$th Kronecker quiver. Then $\mathcal{R}$ is reducible to matrix
normal form of type $\hcf(a,b)$ if and only if $(a\ b)$ is reducible
to matrix normal form for the $t$th Kronecker quiver. When
$\mathcal{R}$ is not reducible to matrix normal form then
$\mathcal{R}$ is a birationally constant family and the general
representation of type $\mathcal{R}$ is isomorphic to
$T_{0}^{a'}\oplus T_{1}^{b'}$ for representations $T_{0}$ and $T_{1}$
such that $\hom(T_{0},T_{1})=\ext(T_{0},T_{1})=\ext(T_{1},T_{0})$ and
$\hom(T_{1},T_{0})=t$ where $a'$ and $b'$ are positive integers such
that $\hcf(a',b')=\hcf(a,b)$.
\end{lemma}
\begin{proof}
If $(a\ b)$ is reducible to matrix normal form for the $t$th Kronecker
quiver we take the family of representations of dimension vector $(a\
b)$ for the $t$th Kronecker quiver and apply the functor $\ \otimes
(S\oplus S')$ considered in lemma \ref{ringelred}. This gives a family
of representations for the quiver with relations $(Q,I)$ that is
birationally representation equivalent to $\mathcal{R}$ since a
general representation of type $\mathcal{R}$ lies in $\E(S,S_{1})$ and
$\mathcal{R}$ is left general with respect to $S$ and therefore shows
that $\mathcal{R}$ is reducible to matrix normal form of type
$\hcf(a,b)$. In the remaining case where $(a\ b)$ is not reducible to
matrix normal form then a general representation of dimension vector
$(a\ b)$ for the $t$th Kronecker quiver is isomorphic to
$R_{0}^{a'}\oplus R_{1}^{b'}$ where
$\hom(R_{0},R_{1})=\ext(R_{0},R_{1})=\ext(R_{1},R_{0})$ and
$\hom(R_{1},R_{0})=t$ and for positive integers $a'$ and $b'$ such
that $\hcf(a',b')=\hcf(a,b)$. Using the facts that
a general representation of type $\mathcal{R}$ lies in $\E(S,S_{1})$
and that $\mathcal{R}$ is left general with respect to $S$, we deduce
that $\mathcal{R}$ is birationally constant and the general
representation of type $\mathcal{R}$ is isomorphic to
$T_{0}^{a'}\oplus T_{1}^{b'}$ where $T_{i}\cong R_{i}\otimes(S\oplus
S')$ and the representations $T_{0}$ and $T_{1}$ and the integers $a'$
and $b'$ satisfy the conditions of the lemma.
\end{proof}

Let $S_{0}$, $S_{1}$ and $S_{2}$ be three representations such that
$\End(S_{i})=k$ and $\Ext(S_{i},S_{i})=0$ for $i=0,1\text{ and }2$,
and $\Hom(S_{i},S_{j})=0=\Ext(S_{i},S_{j})$ for $0\leq i<j\leq 2$. Let
$\mathcal{R}$ be an irreducible family of representations of the
quiver with relations $(Q,I)$. We shall say that $\mathcal{R}$ has a
\definition{two-step Kronecker reduction} of type
$(S_{0},S_{1},S_{2})$ if a general representation of type
$\mathcal{R}$ has an $S_{0}$-socle, $\mathcal{R}$ is left general with
respect to $S_{0}$, and the associated left $S_{0}$-free family has a
Kronecker reduction of type $(S_{1},S_{2})$. Dually, we shall say that
$\mathcal{R}$ has a \definition{two-step Kronecker coreduction} of type
$(S_{0},S_{1},S_{2})$ if a general representation of type
$\mathcal{R}$ has an $S_{2}$-top, $\mathcal{R}$ is right general with
respect to $S_{2}$ and the associated right $S_{2}$-free family has a
Kronecker reduction of type $(S_{0},S_{1})$.  

\begin{lemma}
\label{2step}
Assume that the irreducible family $\mathcal{R}$ of representations of
the quiver with relations $(Q,I)$ over the algebraic variety $X$ has a
two-step Kronecker reduction of type $(S_{0},S_{1},S_{2})$ where
$\hom(S_{0},\mathcal{R})=c$, $\hom(S_{1},\mathcal{R}')=b$ for the
associated left $S_{0}$-free family and $\hom(\mathcal{R}',S_{2})=a$.
Assume that a general representation of type $\mathcal{R}$ has trivial
endomorphism ring. Then $\mathcal{R}$ is reducible to matrix normal
form of type $\hcf(a,b,c)$.
\end{lemma}
\begin{proof}
We shall proceed by induction on $a+b+c$. The associated left
$S_{0}$-free family, $\mathcal{R}'$, has a Kronecker reduction of type
$(S_{1},S_{2})$ to the dimension vector $(a\ b)$ for the $t$th
Kronecker quiver $Q(t)$ for some integer $t$. Now suppose that the
dimension vector $(a\ b)$ is reducible to matrix normal form. Then by
lemma \ref{kr}, $\mathcal{R}'$, is reducible to matrix normal form of
type $\hcf(a,b)$ and so, by lemma \ref{mnfontop}, $\mathcal{R}$ is
reducible to matrix normal form of type $\hcf(a,b,c)$ as required.

If the dimension vector $(a\ b)$ is not reducible to matrix normal
form then, by lemma \ref{kr}, a general representation of type
$\mathcal{R}'$ is isomorphic to $T_{1}^{a'}\oplus T_{2}^{b'}$ where 
$\hom(T_{1},T_{2})=0=\ext(T_{1},T_{2})=\ext(T_{2},T_{1})$ and
$\hom(T_{2},T_{1})=t$ and $a'$ and $b'$ are positive integers such
that $\hcf(a',b')=\hcf(a,b)$. Note that $a'+b'<a+b$ unless $t=0$ in
which case $S_{1}=T_{1}$ and $S_{2}=T_{2}$.

We let $T_{0}=S_{0}$ for notational convenience. A general
representation $R$ of type $\mathcal{R}$ is the middle term of a short
exact sequence
\begin{equation}
\label{dec} 
\ses{T_{0}^{c}}{R}{T_{1}^{a'}\oplus T_{2}^{b'}} 
\end{equation}
Then since $T_{1}$ and $T_{2}$ are representations that lie in
$\E(S_{1},S_{2})$, it follows that
$\hom(T_{0},T_{i})=0=\ext(T_{0},T_{i})$ for $i=1,2$ since
$\hom(T_{0},S_{i})=0=\ext(T_{0},S_{i})$ for $i=1,2$. It also follows
that $\hom(T_{1},T_{0})=0$ from equation \ref{dec} since a general
representation of type $\mathcal{R}$ has trivial endomorphism ring and
$\ext(T_{1},T_{0})\neq 0$ for the same reason.

Since $\hom(T_{i},T_{2})=0$ for $i=0,1$, it follows that
$\hom(R,T_{0})=b'$ and $R$ has a $T_{2}$-top. Let $K$ be the kernel of
the homomorphism from $R$ onto $T_{2}^{b'}$. We have a short exact
sequence 
\begin{equation*}
\ses{T_{0}^{c}}{K}{T_{1}^{a'}}.
\end{equation*}
Therefore the linear map from $\Ext(T_{2}^{b'},T_{0}^{c})$ to
$\Ext(T_{2}^{b'},K)$ is surjective and $\Ext(T_{2}^{b'},T_{0}^{c})$ is
a summand of $\Ext(T_{1}^{a'}\oplus T_{2}^{b'},T_{0}^{c})$.
Therefore, $\mathcal{R}$ is right general with respect to $T_{2}$
because $\mathcal{R}$ is left general with respect to $T_{0}$ and so
an open subvariety of the extensions of $T_{1}^{a'}\oplus T_{2}^{b'}$
on $T_{0}^{c}$ occur in the family $\mathcal{R}$.

Further, the associated right $T_{2}$-free family, $\mathcal{R}''$ has
a Kronecker reduction of type $(T_{0},T_{1})$. The only thing
remaining to check is that the family $\mathcal{R}''$ is left general
with respect to $T_{0}$ but this follows because
$\Ext(T_{1}^{a'},T_{0}^{c})$ is a summand of $\Ext(T_{1}^{a'}\oplus
T_{2}^{b'})$ and the family $\mathcal{R}$ is left general with respect
to $T_{0}$.

Thus we have shown that the family $\mathcal{R}$ has a Kronecker
coreduction of type $(T_{0},T_{1},T_{2})$ and if $a'+b'<a+b$ then
$a'+b'+c<a+b+c$ and we are done by induction. 

If $a'+b'=a+b$, we noted above that $\ext(T_{1},T_{0})\neq 0$ and so
when we perform the same argument again for this coreduction the
numbers will drop this time so again we are done by induction.
\end{proof}

\begin{theorem}
\label{leftgmoduli}
Let $\alpha=(a\ b\ c)$ be a dimension vector for the multiplication
$f$ such that there is a left general representation of dimension
vector $\alpha$ with trivial endomorphism ring. Then the canonical
family on the left general component is reducible to matrix normal
form. Therefore a moduli space of left general representations of this
dimension vector is birational to a suitable number of $h$ by $h$
matrices up to simultaneous conjugacy where $h=\hcf(a,b,c)$.
\end{theorem}
\begin{proof}
Let $\mathcal{R}$ be the canonical family on the left general
component $C$ of representations of dimension vector $\alpha$. Let
$S_{0}= (0\ 0\ k)$, $S_{1}= (0\ k\ 0)$ and $S_{2}=(k\ 0\ 0)$. Then
$\mathcal{R}$ has a two-step Kronecker reduction of type
$(S_{0},S_{1},S_{2})$ and so by lemma \ref{2step}, $\mathcal{R}$ is
reducible to matrix normal form.  
\end{proof}

\section{The main result}
\label{modp2}

By the work of Beilinson \cite{Beilinson}, the derived category of
coherent sheaves on $\mathbb{P}^{2}$ is equivalent to the derived
category of representations of the multiplication $\sigma\colon
k^{3}\otimes k^{3}\rightarrow S^{2}(k^{3})$. We shall identify these
two derived categories using this equivalence. Thus we have a
triangulated category $D$ which has two subcategories
$\Coh(\mathbb{P}^{2})$, the category of coherent sheaves on
$\mathbb{P}^{2}$, and $\Rep(\sigma)$, the category of representations
of the multiplication $\sigma$. We shall make the identification of
the two derived categories in such a way that $\mathcal{O}=(k\ k^{3}\
S^{2}(k^{3}))=P(0)$, $\mathcal{O}(-1)=(0\ k\ k^{3})=P(1)$ and
$\mathcal{O}(-2)=(0\ 0\ k)=P(2)$. The representations $P(0)$, $P(1)$
and $P(2)$ are the indecomposable projective representations of the
multiplication $\sigma$. We shall say that an object of $D$ is a
representation if it lies in $\Rep(\sigma)$ and that it is a sheaf if
it lies in $\Coh(\mathbb{P}^{2})$. Thus we can ask the question
whether a sheaf is a representation and vice versa. One direction is
clear; a sheaf $S$ is a representation if and only if $H^{i}(S(j))=0$
for $i=1,2$ and $j=0,1\text{ and }2$.

\begin{lemma}
\label{reptosheaf}
Let $R=(R(0)\ R(1)\ R(2))$ be a representation of the multiplication
$\sigma$. Let
\begin{align*}
0\rightarrow R(0)\otimes \Lambda^{2}(k^{3})\otimes P(2) &\rightarrow
R(0)\otimes k^{3}\otimes P(1)\oplus R(1)\otimes k^{3}\otimes
P(2)\\
&\rightarrow \oplus_{i=1}^{3}R(i)\otimes P(i)\rightarrow
R\rightarrow 0  
\end{align*}
be a projective resolution of $R$. Then $R$ is a sheaf if and only
if the complex of sheaves
\begin{align*}
0\rightarrow R(0)\otimes \Lambda^{2}(k^{3})\otimes \mathcal{O}(-2)
&\rightarrow R(0)\otimes k^{3}\otimes \mathcal{O}(-1)\oplus R(1)\otimes
k^{3}\otimes \mathcal{O}(-2)\\
&\rightarrow \oplus_{i=1}^{3}R(i)\otimes
\mathcal{O}(-i)\rightarrow 0
\end{align*}
has homology only at the penultimate term.
\end{lemma}
\begin{proof}
Once stated, this result is also clear. In the derived category $D$,
the object $R$ is equivalent to the complex 
\begin{align*}
0\rightarrow R(0)\otimes \Lambda^{2}(k^{3})\otimes P(2) &\rightarrow
R(0)\otimes k^{3}\otimes P(1)\oplus R(1)\otimes k^{3}\otimes
P(2)\\
&\rightarrow \oplus_{i=1}^{3}R(i)\otimes P(i)\rightarrow 0  
\end{align*}
which after our identification is the complex of sheaves
\begin{align*}
0\rightarrow R(0)\otimes \Lambda^{2}(k^{3})\otimes \mathcal{O}(-2)
&\rightarrow R(0)\otimes k^{3}\otimes \mathcal{O}(-1)\oplus R(1)\otimes
k^{3}\otimes \mathcal{O}(-2)\\
&\rightarrow \oplus_{i=1}^{3}R(i)\otimes
\mathcal{O}(-i)\rightarrow 0
\end{align*}
and this itself is a sheaf if and only if it has homology only at the
penultimate term as required.
\end{proof}

This last lemma gives us the following result which allows us to
identify at least birationally families of sheaves and families of
representations.

\begin{lemma}
\label{l:reptosheaf}
Let $\mathcal{R}$ be a family of representations of the multiplication
$\sigma$ over the irreducible algebraic variety $X$. Assume that
there exists a point $p$ such that $\mathcal{R}_{p}$ is a sheaf. Then
there exists an open subvariety $O$ of $X$ such that for all $p\in O$,
$\mathcal{R}_{p}$ is a sheaf. Similarly, if $\mathcal{S}$ is a family
of sheaves on the irreducible algebraic variety $Y$ and there exists a
point $q$ such that $\mathcal{S}_{q}$ is a representation then there
exists an open subvariety $O'$ of $Y$ such that for all $q\in O'$,
$\mathcal{S}_{q}$ is a representation.
\end{lemma}
\begin{proof}
Consider the complex of sheaves on $X\times \mathbb{P}^{2}$ 
\begin{align*}
0\rightarrow \mathcal{R}(0)\otimes S^{2}(k^{3})\otimes \mathcal{O}(-2)
&\rightarrow \mathcal{R}(0)\otimes k^{3}\otimes \mathcal{O}(-1)\oplus
\mathcal{R}(1)\otimes k^{3}\otimes \mathcal{O}(-2)\\
&\rightarrow\oplus_{i=1}^{3}\mathcal{R}(i)\otimes \mathcal{O}(-i)\rightarrow 0.
\end{align*}
Let $Z$ be the support of the homology of this complex except at the
penultimate term. Then $Z$ is closed and so is its image in $X$;
therefore the complement of the image of $Z$ in $X$ is open and it is
the set of points $p$ where $\mathcal{R}_{p}$ is a sheaf.

In the second case, the vanishing of $H^{i}(\mathcal{S}_{q}(j))$ is an
open condition so the result follows.
\end{proof}

We define the \definition{depth} of a vector bundle $E$ over
$\mathbb{P}^{2}$ to be the largest integer $h$ such that $[E]/h$ is in
the Grothendieck group $K_{0}(\Coh(\mathbb{P}^{2}))$ of
$\Coh(\mathbb{P}^{2})$ where $[E]$ is the class of $E$ in this
Grothendieck group. The Grothendieck group of the derived category $D$
coincides with $K_{0}(\Coh(\mathbb{P}^{2}))$ and also with
$K_{0}(\Rep(\sigma))$; hence if $E$ is actually a representation of
dimension vector $(a\ b\ c)$ it follows that the depth of $E$ is
$\hcf(a,b,c)$. 

A sheaf $E$ on $\mathbb{P}^{2}$ is said to have \definition{natural
cohomology} if for all integers $j$ at most one of $H^{i}(E(j))$ is
non-zero for $i=0,1$ and $2$. This definition is important to us since
general vector bundles have this property by \cite{Walter} and
\cite{Hirschowitz} where it comes in the form that slope-semistable
sheaves are prioritary and prioritary sheaves have natural cohomology.
This allows an easy reduction to left general representations.

\begin{theorem}
\label{main}
A moduli space of sheaves on $\mathbb{P}^{2}$ such that the general
sheaf has natural cohomology is birational to a moduli space of left
general representations of the multiplication $\sigma\colon k^{3}\otimes
k^{3}\rightarrow S^{2}(k^{3})$ and hence is birational to a suitable
number of $h$ by $h$ matrices up to simultaneous conjugacy where $h$
is the depth of every sheaf classified by the moduli space. In
particular this holds for a moduli space of vector bundles. 
\end{theorem}
\begin{proof}
After tensoring by a suitable line bundle we may assume for a general
$[S]$ in the moduli space that $S$ is a representation but $S(-1)$ is
not; that is one of $H^{1}(S(-1))$ and $H^{2}(S(-1))$ is
non-zero. Since $S$ has natural cohomology, it follows that
$H^{0}(S(-1))=0$; it also follows that $H^{2}(S(-1))=0$ since by Serre
duality this is dual to $\Hom(S,\mathcal{O}(-2))$ and since $S$ is a
representation any such homomorphism is split surjective.
So we consider the exact complex of sheaves on $\mathbb{P}^{2}$
\begin{equation*}
0\rightarrow \Lambda^{3}(k^{3})\otimes \mathcal{O}(-2)\rightarrow
\Lambda^{2}(k^{3})\otimes \mathcal{O}(-1)\rightarrow k^{3}\otimes
\mathcal{O}\rightarrow \mathcal{O}(1)\rightarrow 0  
\end{equation*}
from which we deduce that for a sheaf $S$ that is also a
representation the groups $H^{i}(S(-1))$ are the homology of the
complex 
\begin{equation*}
0\rightarrow k^{3}{}\dual\otimes H^{0}(S)\rightarrow
\Lambda^{2}(k^{3}){}\dual\otimes H^{0}(S(1))\rightarrow
\Lambda^{3}(k^{3}){}\dual\otimes H^{0}(S(2))\rightarrow 0  
\end{equation*}
and this complex is exact except at the middle term. In particular,
after tensoring with $\Lambda^{3}(k^{3})$ we see that the linear map
from $\Lambda^{2}(k^{3})\otimes H^{0}(S)$ to $k^{3}\otimes
H^{0}(S(1))$ is injective. Let $a=\dim H^{0}(S)$, $b=\dim H^{0}(S(1))$
and $c=\dim H^{0}(S(2))$. Then $S$ considered as a representation has
dimension vector $\alpha=(a\ b\ c)$ and the injectivity of this linear
map shows by lemma \ref{reclg} that $S$ is a representation in the
left general component. Therefore by lemma \ref{l:reptosheaf} our
moduli space of sheaves is birational to an orbit space for
$PGl_{\alpha}$ on the left general component of representations of
dimension vector $\alpha$ but we know that the canonical family on
this component is reducible to matrix normal form by theorem
\ref{leftgmoduli}. Therefore our moduli space is birational to a
suitable number of $h$ by $h$ matrices up to simultaneous conjugacy
where $h=\hcf(a,b,c)$ is the depth of $S$.
\end{proof}

Since the canonical family of left general representations of
dimension vector $(a\ b\ c)$ considered in the proof of this theorem
is reducible to matrix normal form it follows that if we have a moduli
space $M$ of sheaves whose general member has natural cohomology then
there is a family of sheaves in $M$ over an algebraic variety $X$ on
which $PGl_{h}$ acts so that $X$ is $PGl_{h}$ birational to
$M_{h}(k)^{s}$ for some integer $s$ and the morphism from $X$ to $M$
is dominant and is the orbit map.

It is perhaps worth stating the rationality results that follow from
this theorem.

\begin{theorem}
\label{rat}
A moduli space of sheaves on $\mathbb{P}^{2}$ such that the general
sheaf has natural cohomology and depth $n$ is rational when $n=1,2,3
\text{ or }4$. If $4<n$ and $n$ divides $420$ then the moduli space is
stably rational. If $n$ is square-free then the moduli space is
retract rational.
\end{theorem}
\begin{proof}
This follows from the known results on matrices up to simultaneous
conjugacy. A good summary of the known results may be found in
\cite{leB}.  
\end{proof}

\end{document}